
\documentstyle{amsppt}
\magnification=1200
\hoffset=-0.5pc
\nologo
\vsize=57.2truepc
\hsize=38.5truepc

\spaceskip=.5em plus.25em minus.20em

\define\atibottw{4}

\define\singulth{8}
\define\topology{9}
  \define\smooth{10}

 \define\locpois{12}

\define\narashed{18}

\define\sjamaone{20}
\define\sjamlerm{21}

\define\gormacth{23}
\define\geometry{24}
\define\kondsado{25}

\define\raghuboo{26}
\define\ramanato{27}

\noindent
dg-ga/9411007
\newline\noindent
Math. Z., to appear.
\bigskip
\topmatter
\title The singularities of Yang-Mills connections for\\
       bundles on a surface. II. The stratification
\endtitle
\author Johannes Huebschmann{\dag}
\endauthor
\affil
Max Planck Institut f\"ur Mathematik
\\
Gottfried Claren-Str. 26
\\
D-53 225 BONN
\\
huebschm\@mpim-bonn.mpg.de
\endaffil
\date{October 21, 1993}
\enddate
\abstract
{Let $\Sigma$ be a closed surface, $G$ a compact Lie group, not necessarily
connected, with Lie algebra $g$, endowed with an adjoint action invariant
scalar product, let $\xi \colon P \to \Sigma$ be a principal $G$-bundle,
and pick a Riemannian metric and orientation on $\Sigma$ so that the
corresponding Yang-Mills equations are defined. In an earlier paper we
determined the local structure of the moduli space $N(\xi)$ of central
Yang-Mills connections on $\xi$ near an arbitrary point. Here we show that
the decomposition  of $N(\xi)$ into connected components of orbit types of
central Yang-Mills connections is a stratification in the strong  (i.~e.
Whitney) sense; furthermore each stratum,  being a smooth  manifold, inherits
a finite volume symplectic structure from the given data. This complements,
in a way, results of {\smc Atiyah-Bott} in that it will in general decompose
further the critical sets of the corresponding Yang-Mills functional
into smooth pieces.}
\endabstract

\endkeywords
\endsubjclass

\thanks{{\dag} The author carried out this work in the framework
of the VBAC research group of Europroj.}
\endthanks

\endtopmatter
\document
\rightheadtext{The stratification of central Yang-Mills connections}
\leftheadtext{Johannes Huebschmann}

\beginsection Introduction

Let
$\Sigma$ be a closed surface,
$G$ a compact Lie group,
with Lie algebra $g$,
and
$\xi$
a principal $G$-bundle
over $\Sigma$, having a connected total space $P$.
Further, pick
a Riemannian metric on $\Sigma$ and
an {\it orthogonal structure\/}
on $g$,
that is, an adjoint action invariant positive definite
inner product.
These data give rise to a Yang-Mills theory
studied in great detail in \cite\atibottw\ for a connected structure
group. We assume that
solutions of the corresponding Yang-Mills equations exist
--- this will always be the case for a connected
structure group, cf. \cite\atibottw\  --- and we denote
the moduli space of central Yang-Mills connections
by $N(\xi)$.
In
an earlier
paper~\cite\geometry\
to which we refer for background and notation
we
determined the local structure of
$N(\xi)$ near the class of an arbitrary
central Yang-Mills connection.
Extending the approach in~\cite\geometry\
we
establish here the following.

\proclaim{Theorem}
The decomposition
of $N(\xi)$
into
connected components of
orbit types of central Yang-Mills connections
is a Whitney stratification;
furthermore each stratum,
being a smooth  manifold,
inherits
a symplectic structure
from the given data
with finite symplectic volume.
\endproclaim

A more precise statement will be given in (1.2) below,
and the finiteness of the symplectic volumes
will be proved in (1.7) below.
The corresponding decomposition of $N(\xi)$
into pieces will be given in (1.1).
In a way, the Theorem  complements
and extends results of
{\smc Atiyah-Bott}~\cite\atibottw;
in fact,
our stratification
will in general decompose further
each component
of the critical set of
the Yang-Mills functional.
In particular,
for $G=U(n)$, the unitary group,
in the \lq\lq coprime
case\rq\rq, cf. {\smc Atiyah-Bott}~\cite\atibottw, the moduli space
of central Yang-Mills connections
is smooth,
and
our
stratification
then consists of a single piece.
\smallskip
In a follow up paper~\cite\singulth\
we identify the strata
mentioned
in the Theorem
with reductions to suitable
subbundles.
Thereby we {\it cannot\/} avoid running into principal bundles
with {\it non-connected\/} structure groups,
even when the structure group of the bundle
$\xi$ we started with is connected.
This is the reason why the present theory has been set up for
general compact not necessarily connected structure groups.
\smallskip
It is known that the decomposition of
a moduli space of connections
according to orbit types
is a manifold decomposition,
cf. {\smc Kondracki-Sadowski}~\cite\kondsado,
referred to sometimes as a stratification.
The result given in the above Theorem is
really different and
considerably stronger: The space $\Cal N(\xi)$ of central Yang-Mills
connections is {\it not\/} a smooth submanifold of the space of all
connections on $\xi$;
moreover our Theorem  says in particular that the decomposition
into connected components of orbit types is a stratification
in the strong sense, that is, a certain additional
\lq\lq cone condition\rq\rq\ also holds.
\smallskip
The reader is assumed familiar with our paper~\cite\geometry.
Notations and definitions given there
will not be repeated.

\bigskip\noindent
{\bf 1. The stratification}
\smallskip\noindent
We shall use the notion of
{\it stratified space\/}
in the sense of
{\smc Goresky-Mac Pherson}~\cite\gormacth.
For a closed subgroup $K \subseteq G$
let
$\Cal A_{(K)} \subseteq\Cal A(\xi)$
denote the subspace of all connections
$A$ having
stabilizer $Z_A$
whose image in $G$
is conjugate to $K$.
We take as indexing set
$\Cal I$ the set of all possible
stabilizer subgroups
of central Yang-Mills connections $A$ on $\xi$
modulo conjugacy.
For each $(K) \in \Cal I$,
let
$$
\Cal N_{(K)} =\Cal N(\xi) \cap\Cal A_{(K)},\quad
N_{(K)} = \Cal N_{(K)}\big/\Cal G(\xi).
$$
As far as conjugacy classes are concerned
it suffices to take
as representatives of the elements of $\Cal I$
conjugacy classes of subgroups of $G$
rather than
$\Cal G(\xi)$, and we shall henceforth do so.
In particular we shall write
$(K)<G$,
$(K)\in \Cal I$ etc.
The space
$\Cal N(\xi)$ of central Yang-Mills connections
decomposes into a disjoint union
of {\it orbit types\/} $\Cal N_{(K)},\,(K)\in \Cal I$,
and accordingly the moduli space
$N(\xi)$ of central Yang-Mills connections
decomposes into a disjoint union
$$
N(\xi) = \bigcup_{(K) \in \Cal I}N_{(K)}
\tag1.1
$$
of {\it orbit types\/}.
The decomposition we are after is (1.1), with each {\it piece\/}
$N_{(K)}$
decomposed further into its connected components.
We shall not distinguish in notation between (1.1)
and its refinement into connected components.

\proclaim{Theorem 1.2}
The decomposition
of $N(\xi)$
into  connected components
of the pieces of {\rm (1.1)}
is a stratification;
furthermore
each piece $N_{(K)}$,
being a finite union of smooth manifolds,
inherits
a symplectic structure $\sigma_{(K)}$
whose pull back to
$\Cal N_{(K)} \subseteq \Cal A(\xi)$
equals the restriction
to $\Cal N_{(K)}$
of the symplectic form $\sigma$ on $\Cal A(\xi)$.
\endproclaim

\proclaim{Addendum}
The decomposition
of $N(\xi)$
into  connected components
of the pieces of {\rm (1.1)}
is a Whitney stratification.
\endproclaim

\demo{Proof}
Clearly the issue is local.
We shall establish the theorem
by combining
a result
given
in {\smc Sjamaar}~\cite\sjamaone\ (3.5)
and in Section 6
of {\smc Sjamaar-Lerman}~\cite\sjamlerm\
with the local analysis
of  $N(\xi)$
near a point $[A]$ given in Section 2 of our paper~\cite\geometry.
We proceed as follows.
\smallskip
Let $A$ be a central Yang-Mills connection
on $\xi$ fixed henceforth.
Suppose the smooth finite dimensional
symplectic submanifold $\Cal M_A$ of $\Cal A(\xi)$
chosen
as in the proof of
\cite\geometry\ (2.32).
We maintain the notation
$\Cal N_A= \Cal N(\xi)\cap \Cal M_A$
and $N_A= \Cal N_A\big/Z_A$
introduced in Section 2 of \cite\geometry.
The injection of
$\Cal M_A$ into  $\Cal A(\xi)$
induces an injection
of
$N_A$
into
$N(\xi)$
which is in fact a homeomorphism of
$N_A$
onto a neighborhood $U$ (say) of $[A]$ in $N(\xi)$.
This neighborhood $U$, in turn,
inherits a structure of a decomposed
space from the decomposition
(1.1) of $N(\xi)$.
On the other hand,
with respect to the
$Z_A$-structure on
$\Cal N_A$
inherited from that on $\Cal M_A$,
the space
$N_A$
decomposes into connected components of orbit types,
and the
injection
of
$N_A$
into
$N(\xi)$
is decomposition preserving.
\smallskip
Let
$
\vartheta_A
\colon
\Cal M_A
@>>>
z_A^*
$
be the momentum mapping
\cite\geometry\ (2.30) for
the hamiltonian $Z_A$-action
on $\Cal M_A$; it has the value zero at the point $A$.
It is then manifest that
the space
$N_A$
coincides
with the {\it Marsden-Weinstein} reduced space
$\vartheta_A^{-1}(0)\big/Z_A$.
Consequently,
in view of the main result of
{\smc Sjamaar-Lerman}~\cite\sjamlerm,
applied to the momentum mapping $\vartheta_A$
for the {\it compact\/}
group $Z_A$,
the decomposition
of $N_A$
into connected components of orbit types
is a stratification in such a way that
each piece or stratum
inherits a symplectic structure from the given data.
This establishes the theorem since the
injection
of
$N_A$
into
$N(\xi)$
is decomposition preserving. \qed \enddemo
\smallskip
\noindent
{\smc Remark 1.3.}
Let $\Theta_A$
be the momentum mapping
from
$\roman H^1_A(\Sigma,\roman{ad}(\xi))$
to $\roman H^2_A(\Sigma,\roman{ad}(\xi))$
given in \cite\geometry\ (1.2.5),
for the $Z_A$-action on
$\roman H^1_A(\Sigma,\roman{ad}(\xi))$;
it is given by the assignment to
$\eta \in \roman H^1_A(\Sigma,\roman{ad}(\xi))$
of $\frac 12 [\eta,\eta]_A\in\roman H^2_A(\Sigma,\roman{ad}(\xi))$,
where
$[\cdot,\cdot]_A$ refers to the graded bracket
on $\roman H^*_A(\Sigma,\roman{ad}(\xi))$ induced by the data.
By (2.32) of~\cite\geometry,
the reduced space $\roman H_A = \Theta_A^{-1}(0) \big/ Z_A$
is a local model
of a neighborhood of the point $[A]$ of $N(\xi)$.
The theorem can also be established by means of
a direct argument
applied to this local model,
of the kind
used in {\smc Sjamaar}~\cite\sjamaone\ and
{\smc Sjamaar-Lerman}~\cite\sjamlerm\
in the finite dimensional setting.

\demo{Proof of the Addendum}
In Section 3 of
{\smc Sjamaar}~\cite\sjamaone\
and (6.5)
of {\smc Sjamaar-Lerman}~\cite\sjamlerm\
it is proved that a space of the kind
$\roman H_A$
may be embedded into a euclidean space as
a
{\it Whitney\/} stratified set.
This implies that
(1.1)
is a
Whitney
stratification. \qed
\enddemo

\smallskip
The stratification (1.1)
of
$N(\xi)$
has a number of remarkable properties:

\proclaim{Theorem 1.4}
{\rm (1)}
The link of each point of $N(\xi)$, if non-empty,
is connected.
\newline\noindent
{\rm (2)}
There is a unique connected stratum
$N_{(H)}$
which is open and dense in
$N(\xi)$.
\endproclaim

We shall refer to the unique open, connected and dense stratum as
the {\it top\/} stratum of $N(\xi)$, written
$N^{\roman{top}}(\xi)$.

\demo{Proof} Theorem 1.4  is established
in virtually the same way as
the corresponding result
for the reduced space arising from a
momentum mapping for the action of a compact group on
a smooth finite dimensional manifold
due to
{\smc Sjamaar}~\cite\sjamaone\ (3.4.9) and
{\smc Sjamaar-Lerman}~\cite\sjamlerm\ (5.9).
The only difference is that
under the present circumstances
we know {\it a priori\/} that
the reduced space is {\it connected\/}.
We leave the details to the reader. \qed
\enddemo

\smallskip
The following
consequence of (1.4(2))
seems to us to deserve
special mention.

\proclaim{Corollary 1.5}
There is a closed subgroup $T$ of $G$,
unique up to conjugacy,
such that
the space of
gauge equivalence classes of
central Yang-Mills connections $A$
having orbit type $(T)$
coincides with
$N^{\roman{top}}(\xi)$. \qed
\endproclaim

\smallskip
Maintaining terminology
introduced in \cite\geometry, we shall say that
a central Yang-Mills connection $A$
is  {\it non-singular\/} if
its stabilizer $Z_A$ acts trivially on
$\roman H_A^1(\Sigma,\roman{ad}(\xi))$;
the point
$[A]$ of $N(\xi)$
will then be said to be {\it non-singular\/}.
In view of \cite\geometry\ (2.33),
a non-singular central Yang-Mills connection $A$
represents a {\it smooth\/} point of $N(\xi)$
in the sense that its link is a sphere.
It may happen that
the subspace of smooth points of $N(\xi)$
is larger than that of its non-singular ones;
for example this occurs
for $G=\roman {SU}(2)$ over a surface of genus 2,
see \cite\locpois.
However the subspace
of non-singular points
is exactly that
where the symplectic structure
is defined.

\proclaim{Corollary 1.6}
The
space of
gauge equivalence classes of
non-singular central Yang-Mills connections $A$
is non-empty and
coincides with the top stratum
$N^{\roman{top}}(\xi)$.
\endproclaim

\demo{Proof}
Let $A$ be a central Yang-Mills connection.
By
(2.32) of~\cite\geometry,
a neighborhood of
$[A]$
in $N(\xi)$
looks like
a neighborhood of the class of zero in the reduced space
$\roman H_A$.
When $[A]$ lies in the top stratum
of
$N(\xi)$ the class of 0 lies in the top stratum
of
$\roman H_A$, viewed as a stratified symplectic space.
However this is possible only if
$Z_A$ acts trivially on
$\roman H^1_A(\Sigma,\roman{ad}(\xi))$. \qed
\enddemo

\proclaim{Theorem 1.7}
Each stratum of
$N(\xi)$
has finite symplectic volume.
\endproclaim

\demo{Proof}
Let $A$ be a central Yang-Mills connection, let
$\widetilde B_A \subseteq \roman H^1_A(\Sigma,\roman{ad}(\xi))$
be an open relatively compact $Z_A$-invariant ball centered
at the origin,
and let $B_A$ be the reduced space
for the restriction of the momentum mapping
$\Theta_A$ to
$\widetilde B_A$.
As stratified space with symplectic strata,
$B_A$ is a local model
for $N(\xi)$ near the point represented by $A$, cf. (1.3) above.
The argument for (3.9) in
\cite\sjamlerm\
is valid for a space  of the kind
$B_A$ and shows that each stratum
of $B_A$ has finite symplectic volume.
This implies the statement since
$N(\xi)$
may be covered by finitely many open sets having a model
of the kind $B_A$. \qed
\enddemo

\bigskip\noindent
{\bf 2. The relationship with representations}
\medskip\noindent
We pick a base point $Q$ of $\Sigma$.
Let
$
0
@>>>
\bold Z
@>>>
\Gamma
@>>>
\pi
@>>>
1
$
be the universal central extension
of the fundamental group $\pi= \pi_1(\Sigma,Q)$ of $\Sigma$.
It arises from
the standard presentation
$$
\Cal P  = \big\langle x_1,y_1,\dots, x_\ell,y_\ell;
r\big\rangle ,\quad r= \prod_{j=1}^\ell[x_j,y_j],
\tag2.1
$$
of $\pi$
as follows, the number $\ell$ being the genus of $\Sigma$:
Let $F$ be the free group on the generators
and $N$ the normal closure  of $r$ in $F$; then
$\Gamma = F\big / [F,N]$.
Let
$\Gamma_{\bold R}$
be the  group
obtained from $\Gamma$ when its centre
$\bold Z$ is extended to the
additive group $\bold R$ of the reals.
It plays a certain universal role which we now explain:
\smallskip
Write $S^1$ for the circle group, and
let
${\mu\colon  M \to \Sigma}$
be the unique principal $S^1$-bundle
having Chern class 1.
We endow it with a harmonic or Yang-Mills connection
$A_{\Sigma}$
having (normalized) constant curvature
$2\pi i \roman{vol}_{\Sigma}$.
The universal covering projection $\widetilde M \to M$, combined with
$\mu$, yields
the projection map
${
\mu^{\sharp} \colon \widetilde M @>>> \Sigma
}$
of a principal $\Gamma_{\bold R}$-bundle
together with a morphism
of principal bundles
from $\mu^{\sharp}$ to $\mu$;
in particular, it induces a surjective homomorphism
from $\Gamma_{\bold R}$ to $S^1$
inducing an isomorphism of Lie algebras.
Let
$A^{\sharp}_{\Sigma}$
be the lift of
$A_{\Sigma}$
to $\mu^{\sharp}$.
It is a central Yang-Mills connection
on $\mu^{\sharp}$,
with reference to the obvious bi-invariant metric
on
$\Gamma_{\bold R}$.
Let $\Omega \Sigma$ denote the space of piecewise smooth loops in
$\Sigma$, having starting point $Q$;
with the usual composition of loops,
parametrized by intervals of arbitrary length,
$\Omega \Sigma$ is an
{\it associative\/} topological monoid.
\smallskip
We now return to our principal bundle $\xi$.
We pick base points $\widehat Q \in P$
and
$\widetilde Q \in \widetilde M$
over the base point $Q$ of $\Sigma$.
For every connection
$A$ on $\xi$, parallel transport furnishes a continuous
homomorphism
$\tau_{A,\widehat Q}$
of monoids
from $\Omega \Sigma$ to $G$.
In particular,
$A^{\sharp}_{\Sigma}$
induces the homomorphism
$\tau_{A^{\sharp}_{\Sigma},\widetilde Q}$
from
$\Omega \Sigma$
to
$\Gamma_{\bold R}$.
Moreover, cf. \cite\topology\ (2.3),
a Yang-Mills connection $A$ on $\xi$
induces a homomorphism
$\chi_{A,\widehat Q}$ from
$\Gamma_{\bold R}$
to
$G$
so that
$$
\tau_{A,\widehat Q}
=
\chi_{A,\widehat Q}\circ \tau_{A^{\sharp}_{\Sigma},\widetilde Q}
\colon \Omega \Sigma \to G,
$$
whence the assignment
to a central Yang-Mills connection $A$ of
the restriction of the homomorphism
$\chi_{A,\widehat Q}$ to
$\Gamma$ yields a map
from $\Cal N(\xi)$ to
$\roman{Hom}(\Gamma,G)$
the image of which we  denote by
$\roman{Hom}_{\xi}(\Gamma,G)$.
The latter is compact and
closed under conjugation,
and we write
$
\roman{Rep}_{\xi}(\Gamma,G)=
\roman{Hom}_{\xi}(\Gamma,G) \big / G .
$
When $\xi$ is flat
the space $\roman {Hom}_{\xi}(\Gamma,G)$
amounts to the subspace
$\roman {Hom}_{\xi}(\pi,G)$
of $\roman {Hom}(\pi,G)$ corresponding to $\xi$, and
the same kind of remark can be made
for
$\roman {Rep}_{\xi}(\Gamma,G)$.
We maintain the notation $\Cal G^Q(\xi)$
for the group of {\it based\/} gauge transformations.
With these preparations out of the way we recall the following,
cf. Section 2 in our paper~\cite\topology\ for more details.

\proclaim{Theorem 2.3}
The assignment to a central Yang-Mills connection $A$ of $\chi_{A,\widehat Q}$
induces a homeomorphism
from
$\Cal N(\xi)\big/ \Cal G^Q(\xi)$
onto
$\roman {Hom}_{\xi}(\pi,G)$
compatible with the $G$-actions
and hence a homeomorphism
$$
N(\xi)
@>>>
\roman {Rep}_{\xi}(\Gamma,G).
\tag2.3.1
$$
\endproclaim

\smallskip
In view of the identification (2.3.1)
of
$\roman {Rep}_{\xi}(\Gamma,G)$
with
$N(\xi)$,
the stratification (1.1) of
$N(\xi)$
passes to a   decomposition of
$\roman {Rep}_{\xi}(\Gamma,G)$.
The corresponding pieces of
the resulting decomposition of
$\roman {Rep}_{\xi}(\Gamma,G)$
are just as well orbit types,
in view of the following the proof of which
we leave to the reader.

\proclaim{Lemma 2.4}
For every connection $A$ on $\xi$,
the projection map
from $\Cal G(\xi)$ to $G$
(determined by the choice of base point $\widehat Q$)
identifies the stabilizer $Z_A$
of $A$ with the stabilizer
of its class $[A]$ in
$\Cal A(\xi)\big/ \Cal G^Q(\xi)$. \qed
\endproclaim

The Lemma entails that
the bijection (2.3.1)
preserves the decompositions into orbit types on both sides.
\smallskip
Under favorable circumstances,
for example for connected structure group,
the top stratum
of the stratification (1.1) of $N(\xi)$
can be described by representation theory.
Following {\smc Ramanathan}~\cite\ramanato,
we shall say that a representation
$\chi \colon L \to G$
of a group $L$ is {\it irreducible\/}
if the subspace
$g^{L}$
of $L$-invariants
in $g$
under the composite of
$\chi$ with the adjoint representation
of $G$ on $g$
coincides with the
Lie algebra $z$ of the centre $Z$ of $G$.
Now,
for an arbitrary central Yang-Mills connection $A$,
evaluation of a section of
$\roman{ad}(\xi)$ at $Q$ induces an isomorphism from
$\roman H_{A}^0(\Sigma,\roman{ad}(\xi))$
onto the subspace $g^{\Gamma} \subseteq g$
of invariants in   $g$,
with reference to
the induced homomorphism
$\chi_{A,\widehat Q} \in \roman {Hom}_{\xi}(\Gamma,G)$,
and hence $\chi_{A,\widehat Q}$
is manifestly irreducible if and only if
this association
identifies
$\roman H_{A}^0(\Sigma,\roman{ad}(\xi))$
with $z$.
We shall say that
a central Yang-Mills connection $A$
is {\it representation irreducible\/}
if
its stabilizer Lie algebra
$z_A=\roman H_{A}^0(\Sigma,\roman{ad}(\xi))$
is identified
with $z$ in this way.
We have chosen this terminology
in order to avoid conflict
with the common notion of an {\it irreducible\/}
connection.
When $G$ is connected and
the genus of $\Sigma$ is at least 2,
(7.1) and (7.7) of {\smc Ramanathan}~\cite\ramanato\
imply that
there exist
irreducible
representations in $\roman {Hom}_{\xi}(\Gamma,G)$
and hence
representation irreducible
central Yang-Mills connections on $\xi$.

\smallskip
For a representation irreducible
central Yang-Mills connection $A$ on $\xi$,
the stabilizer Lie algebra
$z_A$
amounts to the Lie algebra $z$ of the centre $Z$ of $G$ and the latter
lies in the centre of the Lie algebra $g$;
consequently the action
of the stabilizer $Z_A$
of $A$
on $\roman H_A^1(\Sigma,\roman{ad}(\xi))$
then passes to a representation of the finite group
$\pi_0(Z_A)$ of connected components of $Z_A$ on
$\roman H_A^1(\Sigma,\roman{ad}(\xi))$.

\proclaim{Theorem 2.5}
Suppose that
representation irreducible
central Yang-Mills connections exist.
Then the
top stratum $N^{\roman{top}}$
consists exactly of the
classes $[A]$ of
representation irreducible
central Yang-Mills connections
having the property that
the induced representation
of the finite group
of connected components $\pi_0(Z_A)$ on
$\roman H_A^1(\Sigma,\roman{ad}(\xi))$
is trivial,
whatever representative $A$ of $[A]$.
\endproclaim

\demo{Proof}
Let $A$ be a
representation irreducible
central Yang-Mills connection.
Then the stabilizer
$Z_A$ contains
an isomorphic copy of
the centre $Z$
of $G$ and has Lie algebra $z_A$
isomorphic to the Lie algebra $z$ of $Z$.
Moreover,
cf. Section 2 in~\cite\geometry,
the momentum mapping
$
\Theta_A
$
from
$\roman H_A^1(\Sigma,\roman{ad}(\xi))$
to
$\roman H_A^2(\Sigma,\roman{ad}(\xi))$
is zero.
Hence
a neighborhood of
$[A]$
in $N(\xi)$ looks like
a neighborhood of the image of zero
in the quotient space
$\roman H_A^1(\Sigma,\roman{ad}(\xi))\big/ \pi_0(Z_A)$
modulo the induced action
of the finite group
$\pi_0(Z_A)$ on
$\roman H_A^1(\Sigma,\roman{ad}(\xi))$.
Consequently, cf. (1.6),
the point $[A]$
of
$N(\xi)$
is non-singular if and only
if
the representation
of the finite group
$\pi_0(Z_A)$ on
$\roman H_A^1(\Sigma,\roman{ad}(\xi))$
is trivial. \qed
\enddemo

As observed by {\smc Atiyah-Bott}~\cite\atibottw\
for the case of a connected structure group,
Theorem 2.3 may be used to reduce
the study of central Yang-Mills connections to that of flat
connections.
In fact,
write $Z_e$ for the
connected
component of
the identity of the centre of $G$,
and let $G^{\sharp}= G/Z_e$.
Let
$$
\xi^{\sharp}
\colon
P^{\sharp} =
P/Z_e
@>>>
\Sigma
\tag2.6
$$
be the induced principal $G^{\sharp}$-bundle arising
from dividing out
the group
$Z_e$.
We then have the following.

\proclaim{Theorem 2.7}
The bundle $\xi^{\sharp}$ is flat,
the central Yang-Mills connections on
$\xi^{\sharp}$
are precisely the flat ones,
and the
map from
$N(\xi)$ to
$N(\xi^{\sharp})$
induced by the obvious
morphism
of principal bundles
from $\xi$ to $\xi^{\sharp}$
is in fact
the projection map of a
principal fibre bundle
having compact connected structure group $Z_e^{2\ell} = \roman{Hom}(\pi,Z_e)$.
Moreover this  map is compatible with the
stratifications.
\endproclaim

\smallskip
In particular,
it suffices to study
the stratification of $N(\xi^{\sharp})$.
\smallskip
\demo{Proof}
The first statement has been established in
(3.10) of our paper~\cite\topology.
It is clear that the bundle map is compatible with the
stratifications,
in view of the naturality of the constructions. \qed
\enddemo

\smallskip
\noindent
{\smc Example 2.8.}
In view of
an observation of
{\smc Ramanathan}~\cite\ramanato,
it may well happen
that
the class of
a
representation irreducible
central
Yang-Mills connection $A$
is still a singular point
in the moduli space.
His example is the following one,
cf. (4.1) in ~\cite\ramanato:
Let $H$
be the  subgroup
of $\roman {SO}(n)$
consisting
of diagonal matrices with entries
$\pm 1$;
it is irreducible
in the sense that
the space of $H$-invariants
in the Lie algebra $\roman{so}(n)$
under the adjoint action
is zero.
When the genus $p$ of $\Sigma$ is sufficiently large,
e.~g. $2\ell> 2^n$,
there is a surjective homomorphism
$\chi$ from  $\pi$ to $H$.
The associated
flat connection $A=A_{\chi}$
on the resulting flat principal
$\roman {SO}(n)$-bundle
is representation irreducible
by definition.
The stabilizer $Z_A$ of
$A=A_{\chi}$
coincides with $H$.
Now $\roman H_A^1(\Sigma,\roman{ad}(\xi))$ is
canonically isomorphic to
$\roman H^1(\pi,g_{\chi})$
and the resulting representation
of $H=\pi_0(H)$
on
$\roman H^1(\pi,g_{\chi})$
is non-trivial.
\smallskip
When $\Sigma$ is a 2-torus
there are no irreducible
representations
of $\pi$ in $G$
and hence
no
representation irreducible
central Yang-Mills connections on $\xi$
unless
$G$ is abelian;
for $G$ connected non-abelian,
the top stratum of $N(\xi)$
then corresponds to the conjugacy class $(T)$ of a
maximal torus $T$ in  $G$.

\beginsection 3. Examples

Over the 2-sphere the theory is not interesting since
$N(\xi)$ then consists of a single point.
Until further notice  we assume $G$ connected.
We denote by $H$ the connected component of the identity of the
centre of $G$, and
we write $\ell$ for the genus of $\Sigma$.
\medskip\noindent
{\bf 3.1. Genus 1}
\smallskip\noindent
The fundamental group $\pi = \pi_1(\Sigma)$
is abelian
and there
are no
representation irreducible
central Yang-Mills connections on $\xi$
unless $G$ is abelian, which we exclude henceforth.
Let $T$ be a maximal torus in $G/H$,
and let $W$ be its Weyl group.
It is well known that
the obvious injection
of
$\roman{Hom}(\pi,T)$
into $\roman{Hom}_{\xi^{\sharp}}(\pi,G/H)$
identifies
the space
$(T \times T) \big/ W$,
the $W$-action on
$T \times T$ being the diagonal one,
with the representation space
$\roman{Rep}_{\xi^{\sharp}}(\pi,G/H)$;
see (2.6) for the notation $\xi^{\sharp}$.
Consequently
 $N(\xi)$
is the total space
of a principal $H\times H$-fibre bundle
over
$(T \times T) \big/ W$
as base.
Now
$(T \times T) \big/ W$,
being a $V$-manifold,
is stratified
in the usual way.

\medskip\noindent
{\bf 3.2. Genus $\geq 2$; $G=\roman {SU}(2)$}
\smallskip\noindent
Since the group $G=\roman {SU}(2)$
is simply connected there is only the trivial
$\roman {SU}(2)$-bundle $\xi$ over $\Sigma$.
Consider the space
$\roman{Hom}(\pi,\roman {SU}(2))$
as a subspace of
$G^{2\ell}$ in the usual way, and let
$T$
be the standard circle subgroup inside
$G$; it is a maximal torus.
Consider the space
$\roman{Hom}(\pi,T)$,
viewed as a subspace of $\roman{Hom}(\pi,G)$;
clearly, the former looks like $T^{2\ell}$.
Let $Y$
be the $G$-orbit of
$\roman{Hom}(\pi,T)$
in $\roman{Hom}(\pi,G)$
under the adjoint action.
Then
$\roman{Hom}(\pi,\roman {SU}(2))$
decomposes
into
$\roman{Hom}(\pi,\roman {SU}(2))\setminus Y$
and $Y$.
Each
point in
$\roman{Hom}(\pi,\roman {SU}(2))\setminus Y$
has stabilizer the centre $Z=\{\pm 1\}$,
that is, is an irreducible representation
while
each point in
$\roman{Hom}(\pi,T) \cong T^{2\ell}$
has stabilizer $T$.
Furthermore the inclusion
of
$\roman{Hom}(\pi,T)$ into $Y$
induces a bijection
of
$\roman{Hom}(\pi,T) \big/ W
$
onto
$Y\big/ G$
where $W=\bold Z/2$
refers to the Weyl group of $\roman {SU}(2)$.
Now the non-trivial element $w$ of $W$
acts on
$\roman{Hom}(\pi,T)$, viewed as
$T^{2\ell}$,
by the assignment
to $\left(\zeta_1,\dots,\zeta_{2\ell}\right)\in T^{2\ell}$
of
$\left(\overline\zeta_1,\dots,\overline\zeta_{2\ell}\right)$
where as usual $\overline \zeta$
refers to
the complex conjugate
of $\zeta \in T$.
Hence
the fixed point set $N_{G}$
of
the action of $W$ on
$\roman{Hom}(\pi,T)$
consists of
the homomorphisms $\phi$
having the values
$\pm 1$
on the generators of $\pi$
spelled out in (2.1),
and the $W$-action is free
on
$\roman{Hom}(\pi,T) \setminus N_{G}$.
Thus the resulting stratification looks like
$$
N(\xi) = N_{G} \cup N_{(T)} \cup N_{Z}.
$$
Here
$N_{Z} =(\roman{Hom}(\pi,G)\setminus Y)\big/G$
and
$N_{(T)} =(Y \setminus N_{G})\big/G$.
Moreover the
projection map from
$\roman{Hom}(\pi,G)\setminus Y$
to
$N_{Z}$
is actually a principal $\roman {SO}(3)$-bundle map
whence
$N_{Z}$
is a smooth manifold
of dimension $6\ell-6$.
Furthermore
$N_{(T)}$ is manifestly a connected smooth  manifold
of dimension $2\ell$.
This example is studied in more detail in our papers
\cite\smooth\ and \cite\locpois.

\medskip\noindent
{\bf 3.3. The
trivial $\roman {SO}(3)$-bundle over a surface of genus $\geq 2$}
\smallskip\noindent
Since the group $G=\roman {SO}(3)$
has $\pi_1(G) = \bold Z/2$
there is the trivial $G$-bundle
and a single non-trivial one.
Consider first the
trivial one say $\xi$.
The obvious projection map from
$\roman {SU}(2)$ to $\roman {SO}(3)$
induces a covering projection
from
$\roman{Hom}(\pi,\roman {SU}(2))$
onto
$\roman{Hom}_{\xi}(\pi,\roman {SO}(3))$
having
$\roman{Hom}(\pi,\bold Z/2)$
as its group of deck transformations, that is, the group $(\bold Z/2)^{2\ell}$.
The latter
passes to a covering projection
from
$\roman{Rep}(\pi,\roman {SU}(2))$
onto
$\roman{Rep}_{\xi}(\pi,\roman {SO}(3))$
having still
$(\bold Z/2)^{2\ell}$
as its group of deck transformations.
Under this covering, in turn,
the stratification
$$
N(\widetilde \xi)
=
N_{\widetilde G} \cup
N_{(\widetilde T)} \cup N_{Z}
$$
of
$N(\widetilde \xi)$
for the trivial
$\roman {SU}(2)$-bundle
$\widetilde \xi$
obtained in (3.2) above
(and written $\xi$ there)
passes to the stratification
$$
N(\xi) =
N_{G} \cup
N_{(T)} \cup N_{e};
$$
here
$\widetilde T$ and $T$
refer to the corresponding
circle groups inside
$\widetilde G =\roman {SU}(2)$
and
$G= \roman {SO}(3)$, respectively.

\medskip\noindent
{\bf 3.4. Genus $\geq 2$; $G=\roman {U}(2)$}
\smallskip\noindent
The
connected
 component $H$ of
the identity of
the centre of $U(2)$
coincides with the centre $S^1$
and the quotient
$G/H$ is the group
$\roman {SO}(3)$.
A principal $\roman {U}(2)$-bundle $\xi$ is classified
by its Chern class.
Moreover the corresponding
$\roman {SO}(3)$-bundle
$\xi^{\sharp}$
is trivial or not according as the
Chern class is even or not.
Hence, for even Chern class, in view of
(2.7) above,
$N(\xi)$ is the total space
of a principal $(S^1)^{2\ell}$-bundle
over $N(\xi^{\sharp})$,
and the stratification
of $N(\xi^{\sharp})$
obtained in (3.3)
above induces
the stratification
of
$N(\xi)$.
It is well known,
cf. {\smc Narasimhan-Seshadri}~\cite\narashed,
that, for odd Chern class,
the space $N(\xi)$ is a smooth compact manifold.

\medskip\noindent
{\bf 3.5. The non-trivial $\roman {SO}(3)$-bundle over
a surface of genus $\geq 2$}
\smallskip\noindent
Let $\xi$ be the non-trivial principal
$\roman {SO}(3)$-bundle over $\Sigma$, and
write $\widetilde \xi$ for
a principal
$\roman {U}(2)$-bundle
having odd Chern class.
In view of
(2.7) above,
$N(\widetilde \xi)$ is the total space
of a principal $(S^1)^{2\ell}$-bundle
over $N(\xi)$.
Since
$N(\widetilde \xi)$
is a smooth compact manifold
so is $N(\xi)$.

\medskip\noindent
{\bf 3.6. Genus $\geq 1$; $G=\roman {O}(2)$}
\smallskip\noindent
The group $G=\roman {O}(2)$
is a semi-direct product
$G = S^1 \times_s \bold Z/2$.
Principal $\roman {O}(2)$-bundles
$\xi$ over $\Sigma$
having a connected total space $P$
are classified as follows:
At first,
inspection of its homotopy exact sequence reveals that
$\xi$ determines a unique
homomorphism
$\phi_{\xi}$
from $\pi$ to $\bold Z/2$
which, since $P$ is assumed connected, has to be
{\it non-trivial\/}.
Now, given a non-trivial homomorphism
$\phi$ from
$\pi$ to $\bold Z/2$,
there is always a principal bundle
$\xi$ so that $\phi=\phi_{\xi}$,
and for fixed $\phi$, the principal
$\roman {O}(2)$-bundles
having $\phi_{\xi}=\phi$
are
classified by the cohomology group
$\roman H^2_{\phi}(\Sigma,\pi_1(\roman {SO}(2))$,
where the subscript $\phi$ refers to local coefficients
with respect to the induced
(non-trivial)
$\bold Z/2$-action
on
$\pi_1(\roman {SO}(2)) \cong \bold Z$.
A non-trivial homomorphism
$\phi$ from $\pi$ to $\bold Z/2$
may be realized by a map
from $\Sigma$ to $\bold RP_2$
inducing an isomorphism
of
$\roman H_{\roman{local}}^2(\bold RP_2,\bold Z)$
onto
$\roman H^2_{\phi}(\Sigma,\pi_1(\roman {SO}(2))$.
Under this map we obtain {\it every\/}
principal $\roman {O}(2)$-bundle
on
$\Sigma$ having
$\phi_{\xi} = \phi$
induced from a
principal $\roman {O}(2)$-bundle
on
$\bold RP_2$.
In particular, we see that
the cohomology group
$\roman H^2_{\phi}(\Sigma,\pi_1(\roman {SO}(2))$
classifies these bundles.
\smallskip
We now pick
a
principal $\roman {O}(2)$-bundle
$\xi$
on
$\bold RP_2$
having
$\phi_{\xi} = \phi$
and consider the subspace
$\roman{Hom}_{\xi}(\pi,\roman {O}(2))$
of
$\roman{Hom}(\pi,\roman {O}(2))$
which corresponds to $\xi$.
It admits the following description.
We denote the two elements of
$\bold Z/2$ by $1$ and $-1$, the latter
being the non-trivial element.
Write $G_1=\roman {SO}(2)$,
write
$G_{-1}$ for the other connected component
of $\roman {O}(2)$ and,
with reference to the presentation (2.1),
consider the subspace
$$
G_{\phi(x_1)} \times G_{\phi(y_1)}
\times \cdots \times
G_{\phi(x_\ell)} \times G_{\phi(y_\ell)}
$$
of $G^{2\ell}$.
Then the space
$\roman{Hom}_{\xi}(\pi,\roman {O}(2))$
appears as that of
$2\ell$-tuples
$(u_1,v_1,\dots,u_\ell,v_\ell)$ in
this subspace
having the property
$
\prod [u_j,v_j] = 1.
$
Now, given
a point $\psi$ of
$\roman{Hom}_{\xi}(\pi,\roman {O}(2))$,
the Lie algebra of its stabilizer
amounts to the cohomology group
$\roman H^0(\pi,g_{\psi})$.
However, under the present circumstances,
$g=\bold R$,
and the action of $\pi$ on
$g=\bold R$
factors through the given non-trivial
homomorphism
$\phi$ from $\pi$ to $\bold Z/2$
and hence is non-trivial.
Consequently
$
\roman H^0(\pi,g_{\psi})
=
\roman H^0(\bold Z/2,g_{\psi}) = 0,
$
and the stabilizer of
any point of
$\roman{Hom}_{\xi}(\pi,\roman {O}(2))$,
if non-trivial,
is a finite group.
This implies that
$N(\xi)$ is a symplectic V-manifold
of dimension
equal  to $2\ell-2$.
Thus being a V-manifold,
the space $N(\xi)$, if not smooth, is
stratified
in the usual way.
We do not pursue this further.

\medskip\noindent{\bf 3.7.  Genus $\geq 2$; $G=\roman {O}(3)$}
\smallskip\noindent
The matrix
$-\roman{Id} \in \roman {O}(3)$
is central and has order 2.
Consequently the group $G=\roman {O}(3)$
may be written
as a direct product
$G = \roman {SO}(3) \times \bold Z/2$.
Principal $\roman {O}(3)$-bundles
$\xi$ over  $\Sigma$
having a connected total space $P$
are classified as follows:
As in the previous Subsection,
such a bundle $\xi$ determines a unique
homomorphism
$\phi_{\xi}$ from
$\pi$ to $\bold Z/2$
which, since $P$ is assumed connected, has to be
{\it non-trivial\/}.
Now, given a non-trivial homomorphism
$\phi$ from $\pi$ to $\bold Z/2$,
there is always a principal bundle
$\xi$ so that $\phi=\phi_{\xi}$,
and for fixed $\phi$, the principal
$\roman {O}(3)$-bundles
having $\phi_{\xi}=\phi$
are
classified by the cohomology group
$\roman H^2(\Sigma,\pi_1(\roman {SO}(3))\cong \bold Z/2$.
Thus, given
$\phi$, there are two non-isomorphic principal
$\roman {O}(3)$-bundles
$\xi$ having $\phi_{\xi}=\phi$.
The
direct product decomposition
$\roman {O}(3) = \roman {SO}(3) \times \bold Z/2$
induces a  decomposition
of
$\roman{Hom}(\pi,\roman {O}(3))$
into
the direct product
of
$\roman{Hom}(\pi,\roman {SO}(3))$
and $\roman{Hom}(\pi,\bold Z/2)$.
Hence the stratification of
$N(\xi)$
may be obtained
from the stratification of the corresponding
moduli space for the corresponding
$\roman {SO}(3)$-bundle
described in (3.3) or (3.5) as appropriate,
whatever
principal
$\roman {O}(3)$-bundle $\xi$.
We leave the details to the reader.

\bigskip
\centerline{\smc References}
\medskip\noindent
\widestnumber\key{999}
References 1 -- 22 are given in the first part \cite\geometry.
\medskip

\ref \no  \gormacth
\by M. Goresky and R. MacPherson
\book Stratified Morse Theory
\publ Springer
\publaddr Berlin\linebreak
-Heidelberg-New York-Tokyo
\yr 1988
\endref
\ref \no  \geometry
\by J. Huebschmann
\paper The singularities of Yang-Mills connections
for bundles on a surface. I. The local model
\paperinfo Math. Z. (to appear), dg-ga/9411006
\endref

\ref \no \kondsado
\by W. Kondracki and P. Sadowski
\paper Geometric
structure of the orbit space for the action of automorphisms
on connections
\jour  J. Geom. Phys.
\vol 3
\yr 1986
\pages 421--434
\endref

\ref \no \raghuboo
\by M. S. Raghunatan
\book Discrete subgroups of Lie groups
\publ Springer
\publaddr Berlin-\linebreak
Heidelberg-New York
\yr 1972
\endref

\ref \no \ramanato
\by A. Ramanathan
\paper Stable principal bundles on a compact Riemann surface
\jour Math. Ann.
\vol 213
\yr 1975
\pages  129--152
\endref

\enddocument